\documentclass[12pt]{article}
\usepackage{engl}
\usepackage[dvips]{graphics}

\newenvironment{preuve}{\begin{trivlist}\item[]{\sc Proof}}%
{\nolinebreak $\square$ \end{trivlist}}

\begin{document}

\title{Suppressing nonrevisiting paths}
\author{Fr\'ed\'eric BOSIO}
\maketitle

\begin{abstract}
  In this article, we give, under some hypothesis, a couterexample to the nonrevisiting 
path conjecture and this also might refute conjectures that are known to be equivalent to 
it, especially Hirsch's conjecture.
\end{abstract}

\parbox{13cm}{\parindent 5mm }

\normalsize

\begin{center}
\introduction
\end{center}
  In 1957, W.M. Hirsch asked if every $d$-dimensional polytope with $n$ facets has 
diameter at most $n-d$. This is now referred as Hirsch conjecture. The original question 
was stated also for unbounded polyhedra, but in this case, we know it is false 
\cite{K-W}. Many works have been done in order to prove or refute it and it has been 
established for some special classes of polytopes (\cite{Na},\cite{Kr},\cite{Kl}). We can 
also look at \cite{Zi} for a survey of this theory.

  We know Hirsch conjecture equivalent to the two following conjectures \cite{K-W}: 
\newline
i) The nonrevisiting path conjecture, which states that given two vertices of a 
polytope, there always exists a path joining them whose intersection with any facet is 
connected. \newline
ii) The $d$-step conjecture: Consider a $d$-dimensional simple polytope with $2d$ facets 
and two vertices $x$ and $y$ lying on complementary sets of facets. Then there is a path 
of length $d$ from $x$ to $y$.

  Notice this other formulation of this conjecture: Consider two $(d-1)$-simplices in 
${\Bbb R}^{d-1}$ containing a common point $p$. Then there is a numbering 
$A_1 ,..., A_d $ and $A'_1 ,...,A'_d $ of their respective vertices such that for any 
$i$, $p$ belongs to the convex hull of $A_1 ,..., A_i , A'_{i+1},..., A'_d $.

  On the other hand, many polytopes have been found that are sharp (i.e. meet the bound 
$n-d$) \cite{H-K}, \cite{F-H}, which casts doubt on the validity of the conjecture.

  We show here that, under some hypothesis about combinatorial flips of polytopes, it is 
possible, using a technique that we call path-flipping, to eliminate nonrevisiting paths, 
and that leads to a refutation of the nonrevisiting path conjecture.

\section{Notations and recalls}

  We present here some basic definitions, notations and technics about simple polytopes.

  A polytope means the convex hull of a finite (nonzero) number of points in a real 
affine space. It also can be seen (and preferably for us) as the intersection of some 
closed half-spaces of an affine space, that is bounded and with nonempty interior.

  Here do we only consider simple polytopes (the number of facets containing a given 
vertex is equal to the dimension of the polytope). Also, we are only interested in 
combinatorial polytopes, i.e. isomorphism types of face posets of polytopes. So we will 
call polytope a combinatorial polytope.

\begin{notation}
  The dimesion of a polytope $P$ will be noted $d_P $, its number of facets $n_P $ and 
its number of vertices $n'_P $.
\end{notation}

  We basically indicate vertices of a simple polytope by the set of facets that contain 
it, for example the vertex of a $3$-dimensonal polytope lying on the facets $F$,$F'$ and 
$F''$ ought to be noted $(FF'F'')$.

\begin{definition}
  The two facets that contain an extremity of an edge of a simple polytope an not the 
other one are called the {\em{extremal facets}} of this edge.
\end{definition}

  Notice that a facet containing a vertex of a simple polytope is the extremal facet of 
exactly one edge containing this vertex.

\begin{definition}
  Let $P$ be a polytope. A {\em{path}} in $P$ is a finite sequence $(v_0 ,...,v_l )$ of 
vertices of $P$ such that for any $i, 1 \leq i \leq l$, $v_{i-1}$ and $v_i $ are adjacent 
(i.e. are joined by an edge) in $P$. We also demand that $v_i \neq v_j $ if $i \neq j$.

  The vertex $v_0 $ is called the {\em{origin}} of the path, the vertex $v_l $ its 
{\em{end}} and the integer $l$ its {\em{length}}.

  For $1 \leq i \leq l$, the facet of $P$ (supposed simple) that contains $v_i $ and not 
$v_{i-1}$ is called the $i^{th}$ {\em{arrival facet}} of the path. The facet that 
contains $v_0 $ and not $v_1 $ is called its {\em{start facet}}.
\end{definition}

  If $v$ and $v'$ are vertices of a polytope $P$, their distance is the minimal length of 
a path whose origin is $v$ and whose end is $v'$ (such a path always exist) and we call 
diameter of $P$ the maximal distance between two of its vertices.

\begin{definition}
  A path $(v_0 ,..., v_l )$ in a simple polytope $P$ is called {\em{nonrevisiting}} if, 
for any facet $F$ of $P$, the set of integers $i$ for whose $v_i $ belongs to $P$ is an 
interval of $\Bbb N $.
\end{definition}

\begin{rremarque}
 A path is nonrevisiting if and only if no two of its arrival facets are equal
and none of them contains its origin.
\end{rremarque}

  We now present some manipulations on the simple polytopes. We can refer to \cite{F-H} 
for blending and wedging, or \cite{Ti} for flips.

\subsection{Blendings}
  Consider two simple polytopes $P$ and $Q$ of the same dimension, a vertex $v'_P $ 
(resp. $v'_Q $) of $P$ (resp. of $Q$), and a one-to-one correspondance $\phi $ between 
the facets of $P$ containing $v'_P $ and the facets of $Q$ containing $v'_Q$. We can 
then construct a simple polytope $P \# _{\phi } Q$, or simply $P \# Q$ when no confusion 
is possible, called the blending of $P$ and $Q$ (at $v'_P $ and $v'_Q $ relatively to 
$\phi $) by cutting off small neigbourhoods of the vertices and glueing their 
complements, identifing a facet with its image under $\phi $ (when it has one).

\begin{definition}
  A face of $P$ will be called {\em{stable}} if it doesn't contain $v'_P $ and 
{\em{unstable}} if it does.
\end{definition}

  In the sequel, we won't distinguish stable faces from their natural image under the 
blending.

\begin{proposition}
  We have $n_{P \# Q} = n_P + n_Q - d$ and $n'_{P \# Q} = n'_P + n'_Q - 1$.
\end{proposition}

\begin{definition}
  Consider an edge $e$ in $P$ whose vertices are $v'_P $ and another, noted $v$. Then 
there is exactly one edge $(v,v')$ in $P \# Q $ which contains $v$ and crosses the 
glueing locus. The extremal facet of $(v,v')$ in $P \# Q $ which contains $v'$ is called 
the {\em{virtual extremal facet}} of the edge $e$.

  Consider now a path in $P$ whose origin is $v'_P $. Then the virtual extremal facet of 
its first edge is called the {\em{virtual start facet}} of this path, and the vertex $v'$ 
thereup is called the {\em{virtual origin}} of this path.
\end{definition}

  Though the construction of a blending is symmetric in $P$ and $Q$, our constructions 
won't be. We will call $P$ the left polytope and $Q$ the right one, as well as we call 
things related to $P$ (vertices, stable facets, ...) the left ones and those related to 
$Q$ the right ones.

\subsection{Flips}
  Consider a simple polytope $P$ of dimension at least $3$ and a maximal simplicial face 
$S$ of $P$ (maximal meaning that no greater simplicial face contains it), that is neither 
$P$ nor a facet of $P$. Then we call combinatorial flip of $P$ along $S$ the following 
transformation :

  Let $F_1 ,..., F_k $ the facets of $P$ whose intersection is $S$ and $G_1 ,..., G_l $ 
the others facets that meet $S$. Let $P_S $ the poset obtained by keeping the same set of 
facets as $P$, removing the vertices of $S$ and introducing the vertices that lie on 
every $G_i $ and all but one $F_j $. We then take as face of $P_S $ any intersection of 
facets that contains a vertex. Remark that the vertices that we have introduced are 
pairwise adjacent. We say that $P_S $ is obtained from $P$ by a (combinatorial) flip 
along $S$.

  With our definition, $P$ and $P_S $ have naturally the same set of facets.

N.B. Generally, authors consider that the operations of truncation of a vertex (remark 
such an operation introduces a new facet) or collapsing a simplicial facet onto a vertex 
(if $P$ is not a simplex) also are flips, and sometimes even the operations of appearance 
and disappearance of a simplex. We do not.

\begin{definition}
  Consider a polytope $P$ and a polytope $Q$ obtained from $P$ by a flip along a face 
$E$. Consider a face $F$ of $P$. Its {\em{strict transform}} in $Q$ is the intersection 
(in $Q$) of the facets that contain it (in $P$).

  Besides, we call {\em{new}} a face (and especially a vertex) that is in $Q$ and not in 
$P$ (i.e. the intersection in $P$ of the facets containing this face is empty in $P$).
\end{definition}

  Note that the strict transform of a face $F$ of $P$ is empty if and only if $F$ lies in 
$S$.

\begin{rremarque}
  Let $Q$ a polytope obtained from a polytope $P$ by a flip. Then the new facets of $Q$ 
form a maximal simplicial facet $S'$ and the flip of $Q$ along $S'$ gives back $P$.
\end{rremarque}

\begin{definition}
  An edge (or a simplicial face) of a simple polytope is called {\em{(combinatorially) 
flippable}} if the combinatorial flip of the polytope along this face can be performed an
gives rise to another polytope.
\end{definition}

  The hypothesis that will allow us to construct our counterexample is precisely related 
to flippability of edges.

{\Large{\underline{Principal Hypothesis:}}}
{\it{An edge of a $\geq 3$-dimensional simple polytope is flippable if and only if its 
two extremal facets are disjoint.}}

  The only if part is clear since else the strict transforms of these facets would have a 
nonconnected intersection in the new polytope, the if part is unknown (as far as I know). 
However, there are many cases in which we know it is true. Let's give two representative 
examples:

\begin{proposition}
  The hypothesis is true for edges resulting of a blending (i.e. edges that cross the 
glueing locus).
\end{proposition}

\begin{preuve}
  Consider a polytope $P$ with a vertex $v'_P $, a polytope $Q$ with a vertex $v'_Q $ and 
a bleding $P \# Q$ at these points and an edge $(v,v')$ of the blending that crosses the 
glueing locus, with $v$ in $P$ and $v'$ in $Q$. Consider an hyperplane $H$ whose 
intersection with $P$ is $v'_P $ and perform a projective transformation which sends $H$ 
at infinity. Then $P$ gets an end, and, up to an affine transformation, we can assume 
this end to be the product of a simplex $S$ by a vertical interval $[0 ; +\infty [$. 
Assume furthermore that $v$ is the vertex of $P$ which is sent to the highest place (this 
always can be done by a new affine trasformation that does not change the end). Finally 
truncate the polytope at heigth $\displaystyle{\frac{1}{2}}$.

  Perform the same transformation for $Q$, taking $S \times ]-\infty ; 1[$ as end and 
sending $v'$ to the lowest place, truncating $Q$ at the same height as $P$. The polytope 
obtained by glueing the two pieces is $P \# Q$. Now, from this polytope, let the 
$Q$-piece fall on the $P$-piece. The first two vertices that will met are $v$ and $v'$. 
The polytope we get just after their meeting is obtained by a flip along the desired 
edge. Hence the result.
\end{preuve}

\includegraphics{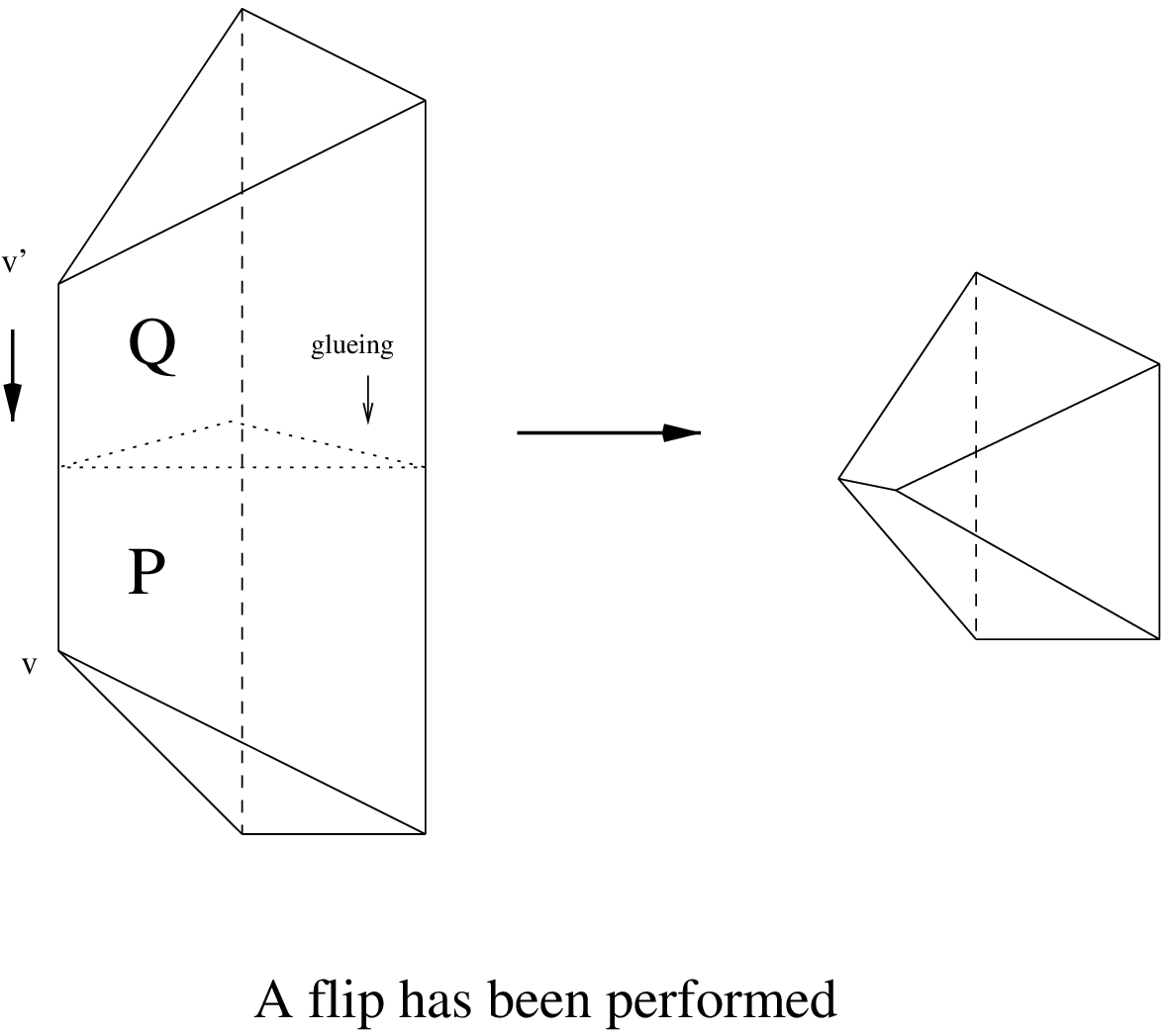}

\begin{proposition}
  The hypothesis is true for $3$-dimensional polytopes.
\end{proposition}

\begin{preuve}
  This results from Steinitz's caracterization of graphs of $3$-polytopes, namely:

\underline{Theorem (Steinitz,1922)}
A graph is the $1$-skeleton of a $3$-dimensional polytope if and only if it is planar and 
$3$-connected (and has $\geq 4$ vertices). $\diamond $

  Obviously, the polytope is simple if and only if the graph is trivalent. Consider the 
graph $G$ of a $3$-dimensional polytope $P$ and an edge $e$ of $G$. Perform now the 
combinatorial flip of $P$ along $e$ and call $G'$ its graph. This can be done except in 
the case, that we will drop, where the edge belongs to a triangular facet of $G$, in 
which case its extremal facets are not disjoint. Thus, we have to show that $G'$ is the 
graph of some $3$-polytope if $e$ has disjoint extremal facets.

  First, we can claim that the new graph $G'$ is planar thanks to the following picture:

\includegraphics{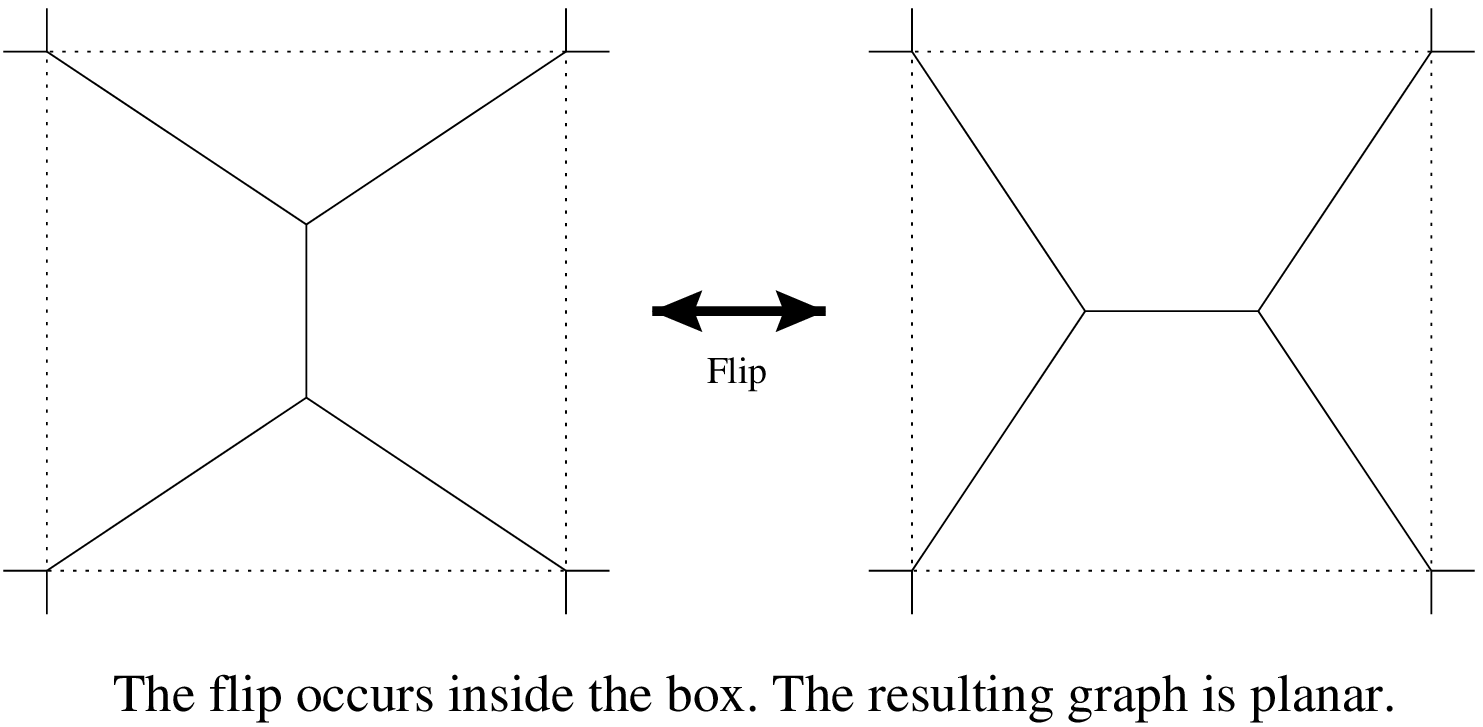}

  If $G'$ is not $3$-connected, consider a pair $\{ v_0 , v_1 \} $ that disconnects him. 
If neither $v_0 $ nor $v_1 $ were one of the two new vertices, then $G$ itself would be 
disconnected by this same pair. If $v_0 $ and $v_1 $ were the two new vertices, $G$ would 
be disconnected by the vertices that have disappeared. Hence we can consider that $v_0 $ 
is a vertex of $G$ and $v_1 $ is new. We can then maintain that the triple $T$ formed by 
$v_0 $ and the two extremities of $e$ disconnects $G$. Now, there is vertex $w$ that is 
adjacent to $v_0 $ and cannot be joined in $G\backslash T$ to any other such vertex. 
Consider finally the two facets of $P$ that contains $v_0 $ and $w$. Each one contains 
one of the two extremities of $e$ else this would contradict what preceeds and none 
contains both as their intersection is the edge $[v_0 ; w]$. So they are the extremal 
facets of $e$ and are not disjoint. The proposition is proved.
\end{preuve}

\subsection{Wedges}

  Let $P$ be a simple polytope with nonempty interior in ${\Bbb R }^n $, $F$ a facet of 
$P$, and $l$ an affine function on ${\Bbb R }^n $ which is nonnegative on $P$ but 
vanishes on $F$. We construct an other polytope $W(P,N)$ by intersecting, in 
${\Bbb R }^n \times \Bbb R $, the $P \times {\Bbb R }_+ $ with the half-space 
$\{ (X,z) \in {\Bbb R }^n \times \Bbb R , z \leq l(X) \} $. This polytope is called the 
wedge over $P$ with foot $F$.

  We generalise this notion by defining, for $k \geq 1$: \\
$W^k (P,F) = \{ (X,z_1 ,..., z_k) \in P \times ({\Bbb R}_+ )^k , 
z_1 \leq z_2 \leq ... \leq z_k \leq l(X) \} $.

\begin{definition}
  The facets that contain $F$ will be called the {\em{large}} facets of $W^k (P,F)$, the 
others well be called the {\em{small}} ones. The large ones will be noted 
$F_i , 0 \leq i \leq k$ and the small ones will be noted like their projections on $P$.
\end{definition}

\begin{proposition}
  A vertex of $W^k (P,F)$ is either a vertex of $F$ or belongs to all but one large 
facet. Moreover, for any vertex $v$ in $P\backslash F$ and large facet $F_i $, there is 
a vertex of $W^k (P,F)$ that projects on $v$ and does not belong to $F_i $.
\end{proposition}

  A vertex of $F$ will be noted by the small facets that contain it, an other vertex will 
be noted by the small facets that contain it, followed by a hat on the large facet that 
does not contain it. For example, consider for $P$ a square, $1,2,1',2'$ its edges. Then, 
in $W^3 (P,1)$ the vertex $(2)$ corresponds to the vertex $(1,2)$ of $P$ and the vertex 
$(1'2\hat{1}_2)$ corresponds to the one which lies on the facets $1_0 $, $1_1 $, $1_3 $, 
$1'$ and $2$.

\begin{corollaire}
  (i) The polytope $W^k (P,F)$ has dimension $d_P + k$, has $n_P + k$ facets, among which 
$k+1$ are large, and has $n'_F + k(n'_P - n'_F)$ vertices; \newline
  (ii) the automorphism group of $W^k (P,F)$ acts transitively on the large facets; 
\newline
  (iii) the wedge over $W^k (P,F)$ with foot a large facet is isomorphic to 
$W^{k+1} (P,F)$.
\end{corollaire}

  The proofs of the proposition an the corollary are left to the reader.

\subsection{The polytope $Q4$}

  There is exactly one combinatorial type of $4$-polytope with $9$ vertices whose 
diameter is equal to $5$. This polytope is called $Q4$ and can be considered as the 
``smallest nontrivial sharp polytope". It has first been discovered by Klee and Walkup 
\cite{K-W}. Here is one of its constructions:

  Consider a $4$-dimensional regular hypercube and choose some point that we consider as 
the ``north pole". The facets containing it are called the upper facets, and we note them 
$1,2,3,4$ whereas their respective opposites, called the lower facets, are noted 
$1',2',3',4'$. Truncate now this hypercube just under the equator. This introduces a 
ninth facet, that we note $N$. Finally perform four flips of edges along the orbit of the 
edge $(1234')$---$(1'234')$ under the natural action of $\Bbb Z _4 $. Now, $Q4$ is the 
resulting polytope. In this polytope, the minimal length of a path joining the vertices 
$(12'34')$ and $(1'23'4)$ is $5$.

N.B.: I like this construction because it illustrates the strength of flipping (and 
especially edge flipping) in suppressing ``small" paths between two points.

  Note the following (easy to verify) property of the polytope $Q4$:

\begin{rremarque}
\label{chemintouch}
  Any nonrevisiting path from $(12'34')$ to $(1'23'4)$ must meet (at least) one of the 
following vertices:

$(2'34'N),(1'2'3N),(1'2'4N),(1'3'4N),(12'4'N),(13'4'N),(23'4'N),(1'23'N)$.
\end{rremarque}

\section{Eliminating nonrevisiting paths}

\subsection{Path-Flips}

  We describe here an efficient way to suppress nonrevisiting paths (as long as the 
hypothesis can be trusted). The idea behind this is that given a nonrevisiting path, it 
is possible to perform a transformation, that we will call path-flip, of the polytope, 
so that one of its vertices now lies on a facet it has already visited. Hence this path 
becomes revisiting.

  Consider a polytope $P$ and a nonrevisiting path $p$ in $P$ such that the start facet 
of $p$ is disjoint from any of its arrival facets. Then we can get a polytope that we 
note $\tilde{P}_p $ in the following manner:

  Assume $p = (v_0 ,..., v_l )$. Then the edge $(v_0 , v_1 )$ is flippable in $P$ as its 
extremal facets, which are the start facet and first arrival facet of $p$, are disjoint. 
In the flipped polytope $\tilde{P}_1 $, $v_2 $ is adjacent to exactly one new vertex 
$v'_1 $. Consider in $\tilde{P}_1 $ the path $p_1 = (v'_1 ,v_2 ,..., v_l )$. It has the 
same start facet as $p$ and its arrival facets are the same (except the first of $p$). 
As the flip of an edge only connects the extremal facets of that edge, the path $p_1 $ in 
$\tilde{P}_1 $ satisfies the same hypothesis as $p$ in $P$. Hence the preceding 
construction can be repeated. This gives finally, after $l$ flips, the polytope 
$\tilde{P}_p $.

  We say that $\tilde{P}_p $ is obtained from $P$ by a path-flip along $p$.

  Assume now that we have several disjoint nonrevisiting paths $p_1 ,..., p_k $ in $P$ 
such that all of them satisfy the formentioned property, no two of them have the same 
start facet and none of them has for start facet an arrival facet of anoter of these 
paths.

  Then we can perform the path-flips altogether, which yields a polytope 
$\tilde{P}_{p_1 ,..., p_k }$.

  To see this by induction, we only have to verify that the strict transforms of 
$p_2 ,...,p_k $ satisfy the same hypothesis in $\tilde{P}_{p_1 }$. This is clear as their 
start and arrival facets have not changed and besides, their start facets cannot having 
been connected with another facet since they are neither the start facet not an arrival 
facet of $p_1 $.
  So $\tilde{P}_{p_1 ,..., p_k }$ exists and is said obtained from $P$ by path-flips 
along the paths $p_1 ,..., p_k $.

\subsection{Vertices without nonrevisiting path}

  Here comes our principal result:

\begin{proposition}
\label{nonrev}
  Consider two polytopes $P$ and $Q$ of the same dimension, points $v_P $, $v'_P $ of $P$ 
and $v_Q $, $v'_Q $ of $Q$, a one-to-one correspondance $\phi $ between the facets of 
$P$ containg $v'_P $ and the facets of $Q$ containg $v'_Q $, and a family of paths $p_i $ 
in $Q$ such that:

\begin{itemize}
\item (i) every $p_i $ has $v'_Q $ as origin;
\item (ii) no $p_i $ contains $v_Q $;
\item (iii) for $i \neq j$, the intersection of $p_i $ and $p_j $, is reduced to $v'_Q $;
\item (iv) every $p_i $ in nonrevisiting;
\item (v) the virtual start facets of the $p_i $ are pairwise different;
\item (vi) in $P$, any nonrevisiting path from $v_P $ to any vertex adjacent to $v'_P $ 
meets every virtual start facet of the $p_i $;
\item (vii) in $Q$, let $v$ be any vertex adjacent to $v'_Q $. Then, if there exists a 
nonrevisiting path from $v$ to $v_Q $ which doesn't meet any $p_i $, then it meets the 
start facet of $(v'_Q , v)$ and $v_P $ lies on the image by $\phi $ of this facet;
\item (viii) for any $i$, $v_Q $ lies on the start facet of $p_i $ and $v_P $ on its 
image by $\phi $.
\end{itemize}

  Then, from $P \# _{\phi } Q$, we can perform path-flips along the paths $p_i $ whose 
origin has been replaced by their virtual origin, and, in the resulting polytope 
$\widetilde{P \# _{\phi } Q}_{p_i }$, there is no nonrevisiting path from $v_P $ to 
$v_Q $.
\end{proposition}

\begin{preuve}
  Let's begin by showing that the path-flips can be performed. We call $p'_i $ the path 
lying in $P \# Q$ where the origin of $p_i $ has been replaced by its vitual origin. For 
any $i$, the start facet of $p'_i $ is a stable left facet and all its arrival facets 
are stable right facets. Hence its start facet is disjoint from all its arrival facets. 
The point $(iv)$ shows that $p'_i $ is nonrevisiting. The paths $p'_i $ are pairwise 
disjoint by $(iii)$. The start facet of $p'_i $ is the left extremal facet of the edge 
containing its origin. They are pairwise different by point $(v)$. Also, every start 
facet of a $p'_i $ is a stable left facet, then cannot be any arrival facet of a $p'_j $ 
since all of them are stable right facets.

  The concomitance of all these properties allows us to claim that the path-flip along 
all the paths $p'_i $ can be constructed, which yields the polytope noted 
$\widetilde{P \# _{\phi } Q}_{p_i }$. Notice also that $v_Q $ is actually one of its 
vertex thanks to point $(ii)$, and that $v_P $ is obviously one also.

  We now prove that every path in $\widetilde{P \# _{\phi } Q}_{p_i }$ whose origin is 
$v_P $ and end is $v_Q $ is revisiting. Consider such a path $p$ and call $v$ its first 
vertex which is not a left vertex. Then:

\underline{First case:} The vertex $v$ is a new vertex. Then, it comes from the flip of 
some path $p_i $. Clearly, $v$ does not lie on the strict transform by the path-flip of 
the facet of $P \# _{\phi } Q$ corresponding to the glueing of the start facet of $p_i $ 
and its image by $\phi $. But by $(viii)$, both $v_P $ and $v_Q $ lie on it. This facet 
is then revisited by $p$.

\underline{Second case:} The vertex $v$ is a right vertex. Then, the first part of $p$, 
from $v_P $ to the vertex immediately before $v$ is in $P$. Then, by $(vi)$, either it is 
revisiting, and then $p$ is, or it meets every virtual left facet, and has left all of 
them when arrived at $v$. Then either $p$ is revisiting or it avoids all the new vertices 
as every such vertex lies on some virtual left facet. By $(vii)$, we can claim that such 
a path meets the strart facet of $(v'_Q ,v')$ where $v'$ is the vertex immeditely 
following the last left vertex ($v'$ might a priori differ from $v$ as the path could 
``return into $P$" after having reached $v$), and even that it meets this facet after 
$v'$. As we also have assumed that $v_P $ lies on this facet, we can maintain that $p$ 
revisits it.

  All in all, every path from $v_P $ to $v_Q $ in $\widetilde{P \# _{\phi } Q}_{p_i }$ is 
revisiting.
\end{preuve}

\section{The counterexample}
  We just now have to find datas that satisfy the conditions of proposition~\ref{nonrev}.
Here are some:

  Consider as $P$ the square of the wedge $W^3 Q4$ and as $Q$ the polytope $W^{10}Q4$. 
Let $v'_P = \left( \begin{array}{c}
1'23'4\hat{N}_0 \\
1'23'4\hat{N}_0 
\end{array} \right) $, 
$v'_Q = (12'34'\hat{N}_0 )$ and $\phi $ the following correspondance of facets: 
$1'_{up, left} \leftrightarrow N_{1,right}, 2_{up, left} \leftrightarrow N_{2,right}, 
3'_{up, left} \leftrightarrow N_{3,right}, 4_{up, left} \leftrightarrow N_{4,right},$ \\
$1'_{down, left} \leftrightarrow N_{5,right}, 2_{down, left} \leftrightarrow N_{6,right}, 
3'_{down, left} \leftrightarrow N_{7,right}, 4_{down, left} \leftrightarrow N_{8,right},$ 
\\ 
$N_{1,up,left} \leftrightarrow 1_{right}, N_{2,up,left} \leftrightarrow 2'_{right},
N_{1,down,left} \leftrightarrow 3_{right}, N_{2,down,left} \leftrightarrow 4'_{right},$ \\
$N_{3,up,left} \leftrightarrow N_{9,right}, N_{3,down,left} \leftrightarrow N_{10,right}$.

  Consider now the following paths in $W^{10}Q4$: \\
$p_1 = (12'34'\hat{N}_0 ) \to (12'34'\hat{N}_9 ) \to (2'34') \to (1'2'3) \to (1'2'4) \to 
(1'3'4)$ and \\
$p_2 = (12'34'\hat{N}_0 ) \to (12'34'\hat{N}_{10}) \to (12'4') \to (13'4') \to 
(23'4') \to (1'23')$. 

  Finally, let $v_P =\left( \begin{array}{c}
12'34'\hat{N}_1 \\
12'34'\hat{N}_1 
\end{array} \right) $
and $v_Q =(1'23'4\hat{N}_0 )$.

\underline{Claim:} These datas satify the hypothesis of propostion~\ref{nonrev}, hence 
this proposition applies, leading to a counterexample to the nonrevisitiong path 
conjecture.

  This is easy to verify. The first three points are obvious. The point (iv) is an 
immediate verification. The virtual start facet of $p_1 $ is $N_{0,up}$, the one of 
$p_2 $ is $N_{0,down}$, hence $(v)$. Besides, both contain $v_P $, hence $(vi)$.

  Consider a vertex $v$ of $Q$ which is adjacent to $v'_Q $. If $v$ has the form 
$(12'34'\hat{N}_k )$, then every nonrevisiting path from $v$ to $v_Q $ meets some $p_i $ 
according to remark~\ref{chemintouch}. If $v$ is on $N$, it is $(12'4')$ (on $p_2 $) or 
$(2'34')$ (on $p_1 $). The last two possibilities are $v = (11'34'\hat{N}_0 )$ or 
$v = (12'33'\hat{N}_0 )$. Consider in this case a nonrevisiting path $p$ from $v$ to 
$v_Q $ that avoids $p_1 $ and $p_2 $. The path $p' = (v'_Q ,p)$ must be revisiting, else 
this would contradict remark~\ref{chemintouch}. The only facet that $p'$ can revisit 
is its start facet and so $p$ meets it. Besides, this facet is $2'$ and its image by 
$\phi $ is $N_{2,up}$ if $v = (11'34'\hat{N}_0 )$ and this facet is $4'$ and its image by 
$\phi $ is $N_{2,down}$ if $v = (12'33'\hat{N}_0 )$. In both cases, $v_P $ lies on this 
facet. So $(vii)$ is verified.

  The start facet of $p_1 $ is $N_9 $, its image by $\phi $ is $N_{3,up}$; the start 
facet of $p_2 $ is $N_{10}$, its image by $\phi $ is $N_{3,down}$. Hence $(viii)$. The 
claim is proved.

\begin{rremarque}
  As $v_P $ and $v_Q $ belong to two common facets, the intersection of these two facets 
is a $12$-dimensional polytope without nonrevisiting path from $v_P $ to $v_Q$.
\end{rremarque}

  Indeed, it is possible that the analyse of the trasformation on this polytope allows us 
to decrease even more the dimension or number of facets of a counterexample.

  Another fact is that it is sometimes believed that the sharpest polytopes are 
neighbourly dual (remark that $Q4$ is such a polytope). The ones we have described are 
not even $2$-neighbourly dual, so it might be possible, following this way, to improve 
the counterexample.

{\footnotesize {Bosio Fr\'ed\'eric \\
Universit\'e de Poitiers \\
UFR Sciences SP2MI \\
D\'epartement de Math\'ematiques \\
UMR CNRS 6086 \\
Teleport 2 \\
Boulevard Marie et Pierre Curie \\
BP 30179 \\
86962 Futuroscope Chasseneuil CEDEX 

e-mail~: bosio@math.univ-poitiers.fr}}

\end{document}